\documentclass[11pt]{article}  
\usepackage{color}
\usepackage{float}
\usepackage{amsmath}
\usepackage{graphicx}
\usepackage{amssymb}
\usepackage{subfig}
\usepackage[latin1]{inputenc}
\usepackage[francais]{babel}
\newcommand{\menge}[2]{\big\{{#1}~\big |~{#2}\big\}}
\newcommand{\RR}{\ensuremath{\mathbb{R}}}
\newcommand{\NN}{\ensuremath{\mathbb N}}
\newcommand{\cart}{\ensuremath{\raisebox{-0.5mm}{\mbox{\LARGE{$\times$}}}}\!}

\begin{document}

\title{\sffamily On the Effectiveness of Projection Methods for 
Convex Feasibility Problems with Linear Inequality Constraints}
\author{Y. Censor,$^{1}$ W. Chen,$^{2}$ P. L. Combettes,$^{3}$
R. Davidi,$^{2}$ and G. T. Herman$^{2}$\\[3mm]
{\small $^{1}$Department of Mathematics, University of Haifa, Mt.
Carmel}
{\small Haifa 31905, Israel}\\
{\small $^{2}$Department of Computer Science, The Graduate Center, City
University of New York}\\ 
{\small 365 Fifth Avenue, New York, NY 10016, USA}\\
{\small $^{3}$UPMC Université Paris 06, Laboratoire Jacques-Louis Lions
-- UMR CNRS 7598}\\ 
{\small 75005 Paris, France}
}

\maketitle

\abstract{The effectiveness of projection methods for solving systems of linear
inequalities is investigated. It is shown that they have a computational
advantage over some alternatives and that this makes them successful
in real-world applications. This is supported by experimental evidence
provided in this paper on problems of various sizes (up to tens of
thousands of unknowns satisfying up to hundreds of thousands of constraints)
and by a discussion of the demonstrated efficacy of projection methods
in numerous scientific publications and commercial patents (dealing
with problems that can have over a billion unknowns and a similar
number of constraints).}

\section{Introduction\label{sec:intro}}

Projection methods were first used to solve systems of linear equations
in Euclidean spaces in the $1930$s \cite{Cimm38,Kacz37} and were
subsequently extended to systems of linear inequalities in \cite{Agmo54,Merz63,Motz54}.
The basic step in these early algorithms consists of a projection
onto an affine subspace or a half-space. Modern projection methods
are much more sophisticated \cite{b96,bb96,bc01,Baus06,Else01,CZ97,Jamo97,Imag97,Hirs05,D01,Svai09,Lopu97}
and they can solve the general \textit{convex feasibility problem}
of finding a point in the intersection of a family of closed convex
sets in a Hilbert space. In such formulations, each set can be specified
in various forms, e.g., as the fixed point set of a nonexpansive operator,
the set of zeros of a maximal monotone operator, the set of solutions
to a convex inequality, or the set of solutions to an equilibrium
problem. Projection methods can have various algorithmic structures
(some of which are particularly suitable for parallel computing) and
they also possess desirable convergence properties and good initial
behavior patterns \cite{bb96,CZ97,c96,Jamo97,Imag97,H09,Pier84}.
The main advantage of projection methods, which makes them successful
in real-world applications, is computational. They commonly have the
ability to handle huge-size problems of dimensions beyond which more
sophisticated methods cease to be efficient or even applicable due
to memory requirements. This is so because the building bricks of
a projection algorithm are the projections onto the given individual
sets, which are assumed to be easy to perform, and because the algorithmic
structure is either sequential or simultaneous, or in-between, as
in the block-iterative projection methods or in the more recently
invented string-averaging projection methods. The number of sets used
simultaneously in each iteration in block-iterative methods and the
number and lengths of strings used in each iteration in string-averaging
methods are variable, which provides great flexibility in matching
the implementation of the algorithm with the parallel architecture
at hand; for block-iterative methods see, e.g., \textcolor{black}{\cite{ac89,Baus06,Butn90,cgg01,Imag97,dhc09,eGG81,Gonz01,Lopu97,Otta88}
and for string-averaging methods see, e.g., \cite{bmr03,bdhk07,ceh01,cs09,ct03,crombez02,pen09,rh03}.}

The convex feasibility formalism is at the core of the modeling of
many problems in various areas of mathematics and the physical sciences;
see \cite{Proc93,c96} and references therein. Over the past four
decades, it has been used to model significant real-world problems
in sensor networks \cite{Hero06}, in radiation therapy treatment
planning \cite{cap88,hc08}, in resolution enhancement \cite{Ceti03},
in wavelet-based denoising \cite{Choi04}, in antenna design \cite{Gu04},
in computerized tomography \cite{H09}, in materials science \cite{kas09},
in watermarking \cite{Lee08}, in data compression \cite{Liew05},
in demosaicking \cite{Lu09}, in magnetic resonance imaging \cite{Sams04},
in holography \cite{Shak08}, in color imaging \cite{Shar00}, in
optics and neural networks \cite{Star98}, in graph matching \cite{Vwik04}
and in adaptive filtering \cite{Yuka06}, to name but a few. In these
-- and numerous other -- problems, projection methods have been used
to solve the underlying convex feasibility problems. 

We focus on the important subclass of convex feasibility problems
in which finitely many sets are given and each of them is specified
by a linear equality or inequality in the Euclidean space $\RR^{N}$.
For such problems, which arise in many important applications \cite{Proc93,H09,hc08},
alternatives to projection methods are available (see, e.g., \cite{A&A,gould}
and the references therein), and it is therefore legitimate to ask
whether projection methods are competitive. 

In this paper we address this question and show that projection methods
are indeed very competitive in the environment of linear inequality
constraints. In Section \ref{sect:comparisons} we discuss their comparative
performance for four different kinds of problems. In Section \ref{sect:published}
we give some examples of their use in real-world applications from
the research and the patent literature. Finally, we present our conclusions.

\section{Comparisons\label{sect:comparisons}}

\subsection{\label{sub:Gould's-2-set-feasibility}Examples of 2-set feasibility
problems }

\label{sec:plc}

In a recent paper \cite{gould}, the author asks in the title: {}``How
good are projection methods for convex feasibility problems?'' and
immediately (in the Abstract) states that:
\begin{quote}
{}``Unfortunately, particularly given the large literature which
might make one think otherwise, numerical tests indicate that in general
none of the variants {[}of projection methods for solving convex feasibility
problems{]} considered are especially effective or competitive with
more sophisticated alternatives.\textquotedblright{}\ 
\end{quote}
As indicated in the Introduction, projection methods have been used
to solve highly nonlinear complex problems involving a very large
number of sets. Therefore, results based on the geometrically simple
2-set problems of \cite{gould} are vastly insufficient to draw general
conclusions. In addition, we show in this subsection that the experiments
reported in \cite{gould} use suboptimal versions of projection methods,
which further questions the justification of the above-quoted general
conclusion as to their effectiveness. 

The numerical experiments provided in \cite{gould} focus exclusively
on the problem of solving a linear system of equations under a box
constraint, namely \begin{equation}
\text{find}\;\; x\in\RR^{N},\;\;\text{such that}\;\;\begin{cases}
Ax=b,\\
x\in\underset{i=1}{\overset{N}{\cart}}\;[c_{i},d_{i}],\end{cases}\label{e:feasgould}\end{equation}
where $A\in\RR^{M\times N}$ ($M\leq N)$ has full rank, $b\in\RR^{M}$,
and the problem is assumed to be feasible. We show that, even in this
basic setting, projection algorithms implemented with standard relaxation
strategies perform much better than indicated by the results in \cite{gould}.

Let us denote by $P_{1}$ and $P_{2}$ the projection operators onto
the closed affine subspace $S_{1}=\menge{x\in\RR^{N}}{Ax=b}$ and
the closed convex set $S_{2}={\cart}_{i=1}^{N}\;[c_{i},d_{i}]$, respectively.
The first operator is defined by \begin{equation}
P_{1}\colon x\mapsto x-A^{\top}\left(AA^{\top}\right){}^{-1}(Ax-b),\label{eq:P1}\end{equation}
where $A^{\top}$ denotes the transpose of $A$. This transformation
can be implemented in various fashions. For instance, in many signal
and image processing problems, the matrix $A$ is block-circulant
and hence diagonalized by the discrete Fourier transform operator,
which leads to a very efficient implementation of $P_{1}$ \cite{Andr77}.
Here, we adopt a QR decomposition approach. Let \begin{equation}
A^{\top}=\left[\begin{array}{cc}
Q_{11} & Q_{12}\end{array}\right]\left[\begin{array}{c}
R{}_{11}\\
0\end{array}\right]\label{eq:QR}\end{equation}
be the QR decomposition of $A^{\top}$, where $R_{11}$ is an $M\times M$
invertible upper triangular matrix \cite{Golu96}. Then (\ref{eq:P1})
yields \begin{equation}
P_{1}\colon x\mapsto x-Q_{11}\left(R_{11}^{\top}\right)^{-1}(Ax-b).\end{equation}
On the other hand, the projection $P_{2}x=(\pi_{i})_{1\leq i\leq N}$
of a vector $x=(x_{i})_{1\leq i\leq N}$ onto $S_{2}$ is obtained
through a simple clipping of its components, i.e., for every $i\in\{1,\ldots,N\}$,
$\pi_{i}=\min\{\max\{x_{i},c_{i}\},d_{i}\}$.

Two standard projection methods to solve \eqref{e:feasgould} are
the alternating projection method \begin{equation}
x^{(0)}\in\RR^{N}\quad\text{and}\quad(\forall n\in\NN)\quad x^{(n+1)}=x^{(n)}+\lambda_{n}(P_{1}P_{2}x^{(n)}-x^{(n)})\label{e:pocs}\end{equation}
and the parallel projection method \begin{equation}
x^{(0)}\in\RR^{N}\quad\text{and}\quad(\forall n\in\NN)\quad x^{(n+1)}=x^{(n)}+\lambda_{n}\bigg(\frac{P_{1}x^{(n)}+P_{2}x^{(n)}}{2}-x^{(n)}\bigg),\label{e:bary}\end{equation}
where $(\lambda_{n})_{n\in\NN}$ is a sequence of strictly positive
relaxation parameters. If $\lambda_{n}\equiv1$ in \eqref{e:pocs},
we obtain the popular \emph{Projection Onto Convex Sets} (POCS) algorithm
\cite{Proc93,Youl82}: \begin{equation}
x^{(0)}\in\RR^{N}\quad\text{and}\quad(\forall n\in\NN)\quad x^{(n+1)}=P_{1}P_{2}x^{(n)}.\label{e:pocs1}\end{equation}
The convergence of any sequence $(x^{(n)})_{n\in\NN}$ thus constructed
to a point in $S_{1}\cap S_{2}$ was established in \cite{Breg65}.
On the other hand, if $\lambda_{n}\equiv1$ in \eqref{e:bary}, we
obtain the \emph{Parallel Projection Method} (PPM): \begin{equation}
x^{(0)}\in\RR^{N}\quad\text{and}\quad(\forall n\in\NN)\quad x^{(n+1)}=\frac{P_{1}x^{(n)}+P_{2}x^{(n)}}{2}.\label{e:ppm}\end{equation}
The convergence of any sequence $(x^{(n)})_{n\in\NN}$ thus constructed
to a point in $S_{1}\cap S_{2}$ was established in \cite{Ausl69},
see also \cite{Aus76}. In \cite{gould}, \eqref{e:pocs} and \eqref{e:bary}
are used, together with variants featuring a construction of $\lambda_{n}$
at iteration $n$ resulting from a line search procedure and without
closed-form expression. However, as the numerical results of \cite{gould}
show, these relaxation schemes do not lead to significantly better
convergence profiles than those obtained with the unrelaxed algorithms
POCS \eqref{e:pocs1} and PPM \eqref{e:ppm}. In addition, nothing
is said regarding the convergence of \eqref{e:pocs} and \eqref{e:bary}
with such relaxation schemes.

The potentially slow convergence of projections methods has long been
recognized \cite{Cott78,Gubi67,Merz63} and remedies have been proposed
to address this problem in the form of adapted relaxation strategies
that guarantee convergence. In the case of \eqref{e:pocs}, it was
shown in \cite{Baus06} that any sequence generated by the \emph{Extrapolated
Alternating Projection Method} (EAPM) \begin{multline}
x^{(0)}\in S_{1}\quad\text{and}\quad(\forall n\in\NN)\quad x^{(n+1)}=x^{(n)}+\rho K_{n}(P_{1}P_{2}x^{(n)}-x^{(n)}),\\
\quad\text{where}\quad0<\rho<2\quad\text{and}\quad K_{n}=\begin{cases}
{\displaystyle \frac{\|P_{2}x^{(n)}-x^{(n)}\|^{2}}{\|P_{1}P_{2}x^{(n)}-x^{(n)}\|^{2}},} & \text{if}\;\; x^{(n)}\notin S_{2},\\
1, & \text{if}\;\; x^{(n)}\in S_{2},\end{cases}\label{eq:e:eapm}\end{multline}
produces a fast algorithm that converges to a solution to \eqref{e:feasgould}.
This type of extrapolation scheme, which exploits the fact that $S_{1}$
is an affine subspace, actually goes back to the classical work of
\cite{Gubi67}. It has been further investigated in \cite{Baus03,Imag97}
and has been extended recently to a general block-iterative scheme
in \cite{Baus06}. Acceleration methods have also been devised for
the parallel algorithm \eqref{e:bary}. Thus, the convergence of the
sequence produced by the \emph{Extrapolated Parallel Projection Method}
(EPPM) \begin{multline}
x^{(0)}\in\RR^{N}\quad\text{and}\quad(\forall n\in\NN)\quad x^{(n+1)}=x^{(n)}+\chi L_{n}\bigg(\frac{P_{1}x^{(n)}+P_{2}x^{(n)}}{2}-x^{(n)}\bigg),\\[2mm]
\quad\text{where}\quad0<\chi<2\quad\text{and}\hskip64mm~\\
L_{n}=\begin{cases}
2{\displaystyle {\frac{\|P_{1}x^{(n)}-x^{(n)}\|^{2}+\|P_{2}x^{(n)}-x^{(n)}\|^{2}}{\|P_{1}x^{(n)}+P_{2}x^{(n)}-2x^{(n)}\|^{2}}},} & \text{if}\;\; x^{(n)}\notin S_{1}\cap S_{2},\\
1, & \text{if}\;\; x^{(n)}\in S_{1}\cap S_{2},\end{cases}\label{eq:e:eppm}\end{multline}
to a solution of \eqref{e:feasgould} was established in \cite{Jamo97}.
This type of parallel extrapolated method goes back to\cite{Merz63}
and \cite{Pier76} , and it has been refined or generalized in several
places \cite{Imag97,Lopu97,Otta88}. In particular, it has been shown
in numerical experiments to be much faster than unrelaxed projection
algorithms in various types of problems ranging from numerical PDEs
to image processing \cite{Imag97,Gonz01,Pier76,Pier84}.

\begin{figure}[t]
\centering{}$\qquad\qquad$\includegraphics[bb=50bp 150bp 612bp 592bp,scale=0.47]{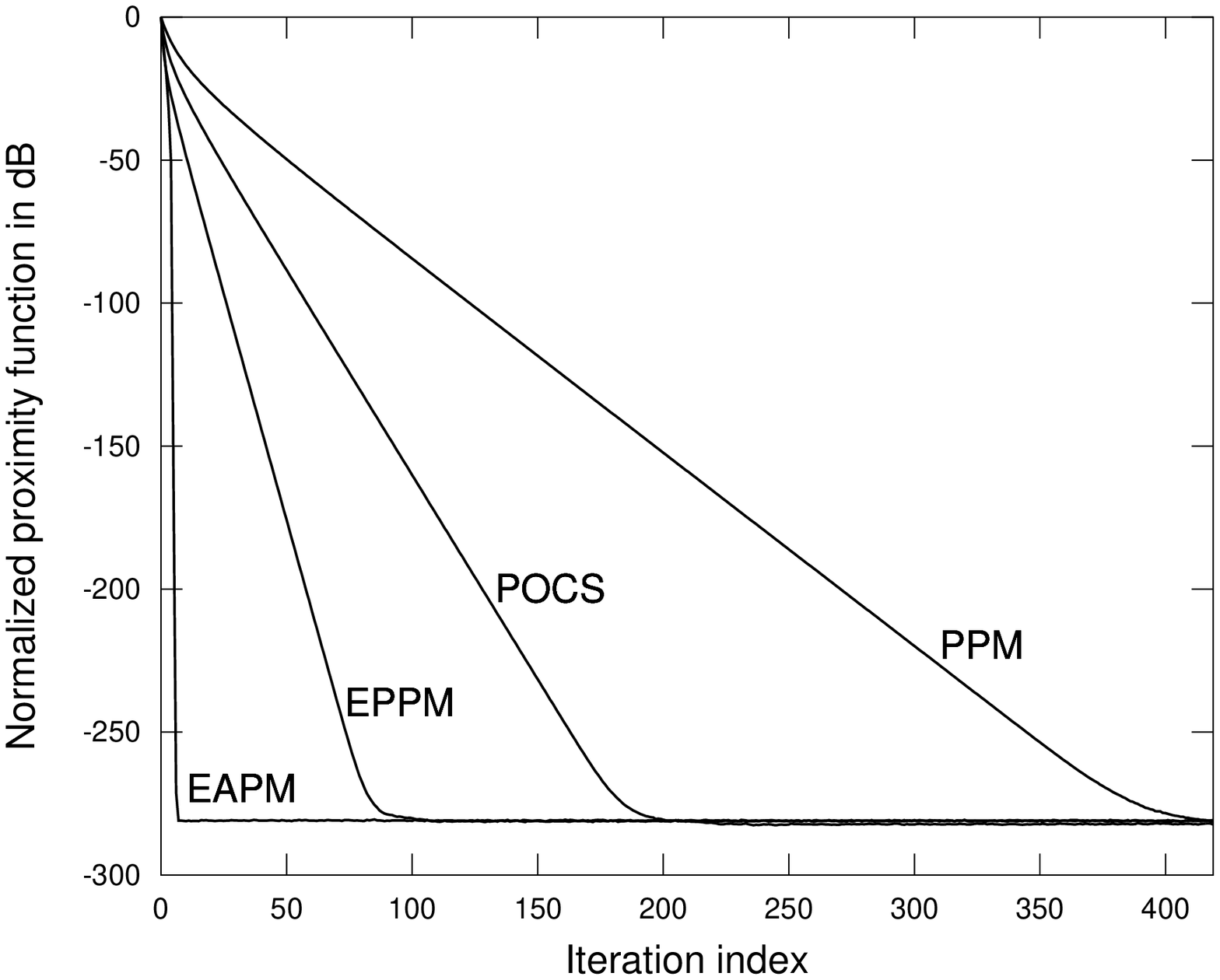}
\caption{\label{fig:.combettes_1}Average performance of the algorithms when
$M\times N=600\times1000$.}

\end{figure}

\begin{figure}
$\qquad\qquad\qquad$$\quad$\includegraphics[bb=50bp 150bp 612bp 592bp,scale=0.47]{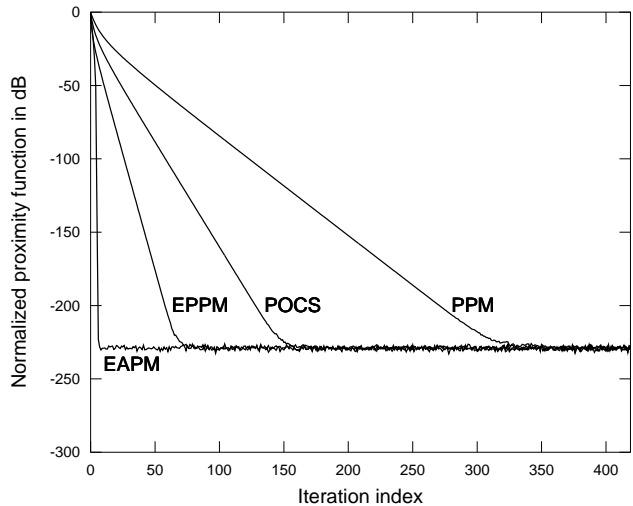}

\caption{\label{fig:.combettes_1-1}Average performance of the algorithms when
$M\times N=600\times1000$ and condition numbers are around $3\times10^{4}$.}

\end{figure}

In Figure~\ref{fig:.combettes_1}, we compare the numerical performance
of POCS \eqref{e:pocs1}, PPM \eqref{e:ppm}, EAPM \eqref{eq:e:eapm},
and EPPM \eqref{eq:e:eppm} for problems of size $M\times N=600\times1000$.
As in \cite{Baus06,Imag97,c96}, the performance of the algorithms
is measured by the decibel (dB) values of the normalized proximity
function, which is evaluated at the $n$th iterate $x^{(n)}$ by \begin{equation}
10\log_{10}\left(\frac{\|P_{1}x^{(n)}-x^{(n)}\|^{2}+\|P_{2}x^{(n)}-x^{(n)}\|^{2}}{\|P_{1}x^{(0)}-x^{(0)}\|^{2}+\|P_{2}x^{(0)}-x^{(0)}\|^{2}}\right).\label{e:proximity}\end{equation}
This comparison is relevant because the computational load of each
iteration resides essentially in the computation of the projection
onto $S_{1}$ and it is therefore roughly the same for all four algorithms.
The results are averaged over 20 runs of the algorithms initialized
with $x^{(0)}=P_{1}0$ and $\rho=\chi=1.9$. In each run a matrix
$A\in[-0.5,0.5]^{M\times N}$ and a vector $x\in[0,1]^{N}$ are randomly
generated. The vector $b=Ax$ is then constructed so as to obtain
a feasible problem using $c_{i}\equiv0$ and $d_{i}\equiv1$ in \eqref{e:feasgould}.
As in \cite{gould} and many other studies, we observe that POCS is
faster than PPM. However, EPPM is faster than POCS and EAPM is clearly
the best method: on the average, it is about $60$ times faster than
PPM, $30$ times faster than POCS, and it achieves full convergence
in just $7$ iterations. In addition, convergence to a feasible solution
is guaranteed by the theory and the expression of the extrapolation
parameter $K_{n}$ in (\ref{eq:e:eapm}) is explicit and it requires
no additional computation. It is argued in Section 5 of \cite{gould}
that ``there is a significant difference between random and real-life
problems (similar observations have been made for linear equations,
where random problems tend to be well-conditioned {[}Reference{]},
and thus often easier to solve than those from applications).'' Let
us observe that random matrices do show up in many real-life problems,
see \cite{Sign89,Zhan05} and the references therein. In addition,
as shown in Figure \ref{fig:.combettes_1-1}, the qualitative behavior
of the algorithms in the presence of poor conditioning is quite comparable
to that observed in Figure \ref{fig:.combettes_1} (for the experiments
of Figure \ref{fig:.combettes_1-1}, the condition numbers vary from
$3\times10^{4}$ to $3.5\times10^{4}$). 

\begin{figure}[H]
\centering{}$\qquad\qquad$$\quad$\includegraphics[bb=50bp 150bp 612bp 592bp,scale=0.47]{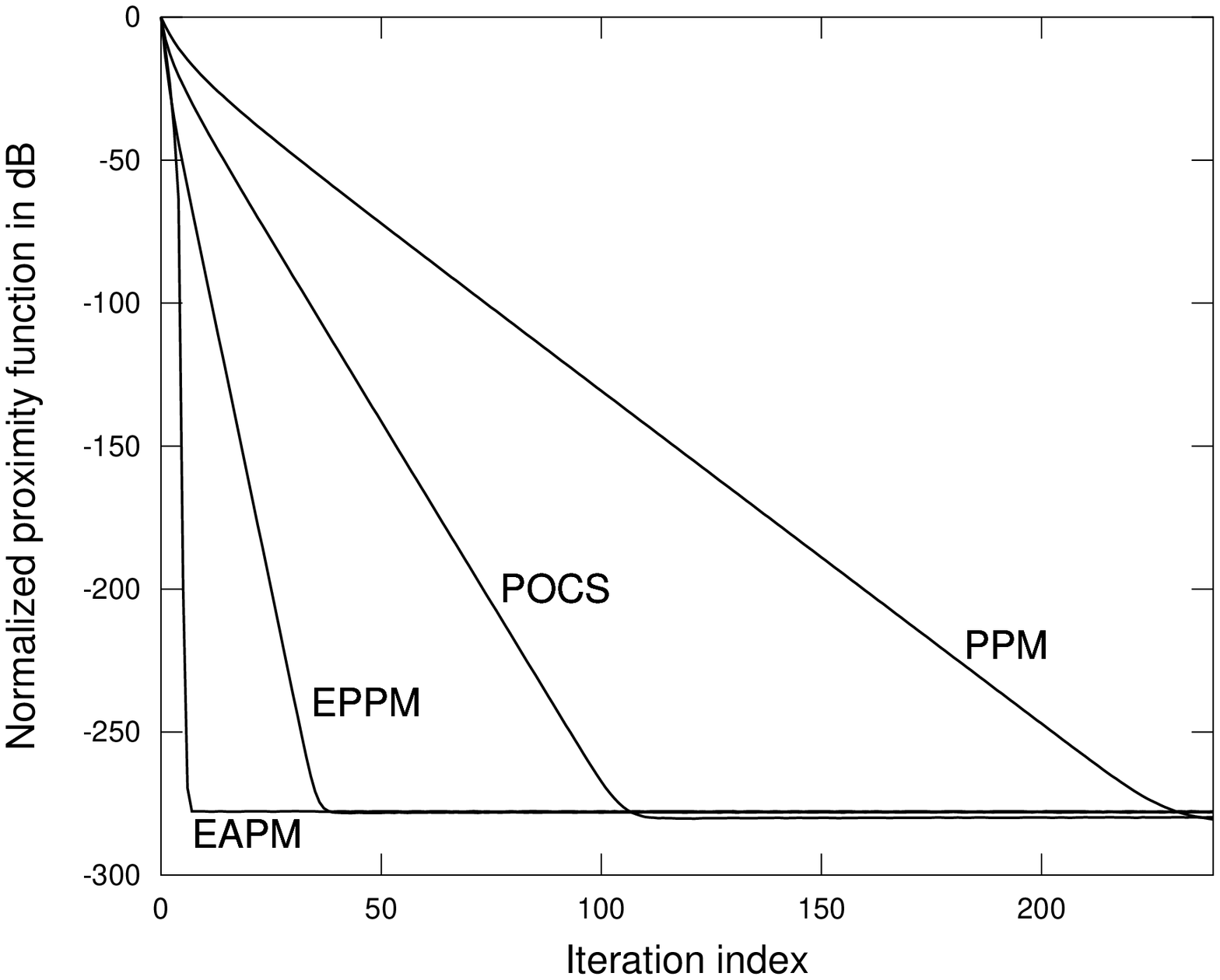}\caption{\label{fig:Combettes_2}Average performance of the algorithms when
$M\times N=3000\times7000$.}

\end{figure}

We have consistently observed this type of performance for problems
of various sizes. For instance, we report in Figure~\ref{fig:Combettes_2}
on the same experiment as above on problems of size $M\times N=3000\times7000$.
Here EAPM is about $45$ times faster than PPM, $22$ times faster
than POCS, and full convergence is achieved in just $5$ iterations.

These experiments indicate that the results in \cite{gould} on the
speed of convergence of POCS and PPM (and the variants proposed there
featuring modest speed-up factors and lacking a formal convergence
analysis) correspond to a suboptimal implementation of projection
methods and are not representative of their performance, since drastic
improvements can be achieved by appropriate relaxations.

\subsection{\label{sub:An-example-from}An example from image representation}

The problem with the largest number of unknowns in the Netlib/CUTEr
LP problem set used in \cite{gould} has $M\times N=6,330\times22,275$
and (according to the on-line attachment to \cite{gould}), for that
problem, all methods discussed in \cite{gould} need 42 seconds or
more to reach the stopping tolerance on a $3.06$ GHz Dell Precision
$650$ workstation. We found among the problems from applications
that we have been investigating one that is over an order of magnitude
larger and for which the projection algorithm recommended in \cite{Chen&Herm}
required only $25$ seconds on the average on an Intel Xeon $1.7$
GHz processor, $1$ Gbyte memory 32 bit workstation using the SNARK09
programming system \cite{snark09}. We now give a brief description
of this problem.

A $J\times J$ digitized image is one that is subdivided into $J^{2}$
square-shaped pixels within each of which the image value is uniform.
Sometimes alternative representations of an image are superior. For
example, in computerized tomography \cite{H09}, we use the blob basis
functions advocated by Lewitt \cite{Lew90} in some series expansion
methods to reduce artifacts in the reconstruction. Such a reduction
is due to the fact that blob basis functions are smoother than pixel
basis functions.

The contribution to the image value at the center of any of the $M=J^{2}$
pixels by any of the $N$ blob basis functions is known from the geometry
of the representations. If we are given a pixel image to start with
and would like to find a good blob representation for it, the task
is to find the weights $x$ to be given to the blobs so that their
combined contributions approximate the pixel values. In mathematical
terms, this problem can be formulated as 

\begin{equation}
\text{find}\;\; x\in\RR^{N}\;\;\text{such that}\;\; c\leq Ax\leq d,\label{eq:blob-representation}\end{equation}
where the bounds $c$ and $d$ have to be tight to ensure a good approximation
of the pixel image by the blob image. (The entries in the matrix $A$
are the values of the various blobs at the centers of the various
pixels.) 

In the experiments reported in \cite{Chen&Herm} $M\times N=59,\!049\times51,\!152$.
The algorithm that was found most efficacious among those tried is
the projection method called CART3$^{\text{++}}$: the average (over
$40$ instances of the problem) time required by CART3$^{\text{++}}$
to find a solution to \eqref{eq:blob-representation} was less than
$25$ seconds.

The algorithm CART3$^{\text{++}}$ belongs to a large family of projection
methods that are usually referred to as \emph{algebraic reconstruction
techniques }(ART). These were first introduced to the tomographic
image reconstruction literature in \cite{gbh70}; for a recent discussion,
see \cite[Chapter 11]{H09}. CART3$^{\text{++}}$, just like the closely
related ART3+ that is used to solve the problem discussed in the next
subsection, has an interesting mathematical property: provided that
the set of feasible vectors satisfying the inequalities in (\ref{eq:blob-representation})
has a nonempty interior, both CART3$^{\text{++}}$and ART3+ will find
a feasible solution in a finite number of iterations \cite{Chen&Herm}.
This is achieved by appropriately controlling the sequence of relaxation
parameters associated with the individual projections.

\subsection{An example from intensity-modulated radiation therapy planning}

The goal of intensity-modulated radiation therapy is to deliver sufficient
doses to tumors to kill them, but without causing irreparable damage
to critical organs. This requirement can be formulated as a linear
feasibility problem of the kind shown in (\ref{eq:blob-representation}).
The interpretation in this application is that each component of $x$
is a to-be-determined strength of radiation to be delivered to the
patient in $N$ separate beamlets, the components of $Ax$ are the
resulting doses at $M$ points in the patient's body, and $c$ and
$d$ are provided by the radiation oncologist as the desired limits
on these doses. Two of the authors of the present paper (W. Chen and
G.T. Herman) have been working in this area with D. Craft, T.M. Madden,
K. Zhang and H.M. Kooy of the Department of Radiation Oncology, Massachusetts
General Hospital and Harvard Medical School, and what follows in this
subsection is an outcome of this collaboration.

In the clinical case that we use as an example we have $M\times N=302,\!491\times13,\!734$.
The number of nonzero elements in $A$ is $62,\!226,\!127$, which
is less than 1.5\% of the total number of entries of $A$, an important
consideration for the efficacy of projection methods for solving the
problem. There is an additional technical consideration: since it
is impossible to deliver negative radiation, each component of $x$
has to be nonnegative, which results in an additional $13,\!734$
inequality constraints. As mentioned at the end of the last subsection,
we use ART3+ \cite{hc08} to solve this feasibility problem.

In clinical applications, it is considered desirable to find multiple
feasible points, each of which is optimal according to its own criterion.
A typical optimization task is {}``find a feasible point that results
in the smallest total dose delivered to the liver.'' The associated
functional is a linear one: it is the sum of those components of $Ax$
that are associated with points in the liver. Recognizing the speed
by which ART3+ finds a feasible point, we propose to apply it repeatedly,
to solve the linear optimization problem

\begin{equation}
\text{find}\;\; x\in\RR^{N}\;\;\mbox{that minimizes \ensuremath{a^{\top}x\,\,}subject to }c\leq Ax\leq d.\label{eq:LP}\end{equation}
Our method solves this problem by turning the objective function into
an additional constraint and solving

\begin{equation}
\text{find}\;\; x\in\RR^{N}\;\;\text{such that}\;\; c\leq Ax\leq d\mbox{ and \ensuremath{a^{\top}x\le\rho}}\label{eq:extra constraint}\end{equation}
using ART3+. By reducing $\rho$ using a bisection search until we
obtain (within a prespecified tolerance) the lowest value possible
for it, we get a good approximation to a solution of (\ref{eq:LP}).
This whole process is called ART3+O.

\begin{figure}[t]
\begin{centering}
\includegraphics[scale=0.9]{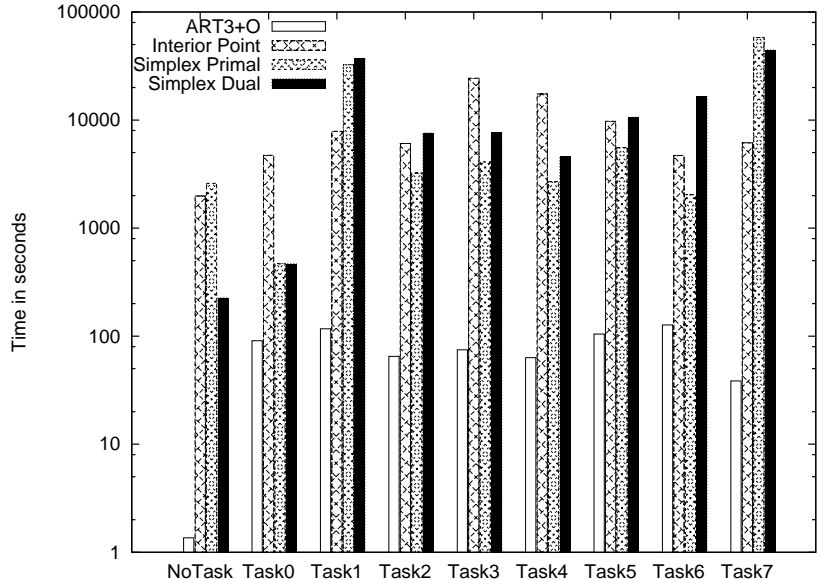}
\par\end{centering}

\caption{\label{timing}Timings of the four methods for the feasibility run
and the eight optimization tasks.}

\end{figure}

The task of minimizing a linear functional subject to linear inequality
constraints is the well-known \emph{Linear Programming }(LP) problem
and several software packages are available for solving it, see, e.g.,
\cite{A&A}. To compare the efficiency of our proposed procedure with
currently popular standard approaches, we applied them to the problem
(\ref{eq:LP}) for a patient with pancreatic cancer. We used all methods
to find just a feasible point (No Task) and also for eight different
LP tasks representing various linear optimization criteria. The three
algorithms with which we compared ART3+O were the self-dual interior
point optimizer, the primal simplex optimizer and the dual simplex
optimizer in the commercial software package MOSEK version 5. The
results are reported in Figure \ref{timing}. Typically, for each
task, ART3+O used about one to two minutes and the MOSEK algorithms
needed one to several hours on an Intel Xeon 2.66 GHz processor, 16
Gbyte memory, 64 bit workstation. It is also noteworthy that the memory
requirements of the MOSEK algorithms were at least twelve times as
large as that of ART3+O. 

\subsection{Examples from computerized tomography}

Computerized tomography is the problem of recovering an image from
its measured (and hence not strictly accurate) integrals along $M$
lines \cite{H09}. If we assume that the recovered image will be represented
as a linear combination of $N$ basis functions (see Subsection \ref{sub:An-example-from}),
then the task is to find the vector $x$ the components of which are
the weights to be given to the basis functions. Due to the linearity
of integration and based on the knowledge of the basis functions,
we can produce an $M\times N$ matrix $A$ such that $Ax$ is approximately
the vector $b$ of measurements. Since it is not likely that there
is an $x$ such that $Ax=b$, it is reasonable to aim instead at finding
an $x$ that minimizes

\begin{equation}
\sigma^{2}\left\Vert b-Ax\right\Vert ^{2}+\left\Vert x\right\Vert ^{2},\label{eq:minimize}\end{equation}
where $\sigma\in\RR$ indicates our confidence in our measurements.
As explained in Section 11.3 of \cite{H09}, this sought-after $x$
is in fact the $x$ part of the minimum norm solution of the consistent
system of equations

\begin{equation}
\left[U\;\sigma A\right]\left[{u\atop x}\right]=\sigma b,\label{eq:normal}\end{equation}
where $U$ is the $M\times M$ identity matrix. In the same section
there is a derivation of a variant of ART that converges to the sought-after
$x$, given by:

\begin{equation}
\begin{array}{l}
u^{(0)}\,\,\textnormal{is the}\,\, M\textnormal{-dimensional zero vector,}\\
x^{(0)}\,\,\textnormal{is the}\,\, N\textnormal{-dimensional zero vector},\\
u^{(n+1)}=u^{(n)}+\gamma_{n}e_{j_{n}},\\
x^{(n+1)}=x^{(n)}+\sigma\gamma_{n}a_{j_{n}},\end{array}\label{U^(0)   , X^(0) is mu_x, U^(k+1),  X^(K+1)  Eq}\end{equation}
with\begin{equation}
\gamma_{n}=\lambda\frac{\sigma\left(b_{j_{n}}-a_{j_{n}}^{\top}x^{(n)}\right)-u_{j_{n}}^{(n)}}{1+\sigma^{2}\left\Vert a_{j_{n}}\right\Vert ^{2}},\label{c^(k) relate to U^(0)   , X^(0) is mu_x, U^(k+1),  X^(K+1)  Eq}\end{equation}
where, for $n\in\NN$, $j_{n}=(n\mbox{ mod }M)+1$, for $1\leq j\leq M$,
$e_{j}$ is the $M$-dimensional vector whose $j$th component is
$1$ and whose other components are $0$, $a_{j}^{\top}$ is the $j$th
row of $A$ and $b_{j}$ is the $j$th component of $b$, and $0<\lambda<2$.
Recognizing that in one iterative step only one row of the matrix
is needed and that in computerized tomography most entries of each
row are zero, we see that an iterative step can be carried out very
rapidly, provided that we have access to the locations and the values
of the nonzero entries. If the memory of the computer is large enough,
this can be accommodated by storing $A$ in a row-by-row sparse representation,
otherwise the locations and values of the nonzero entries can be generated
within each iterative step by some rapid mechanism, such as the \emph{digital
difference analyze}r explained, e.g., in Section $4.6$ of \cite{H09}.

\begin{figure}
\subfloat[Plot of the objective function (\ref{eq:minimize}).]{\includegraphics[scale=0.9]{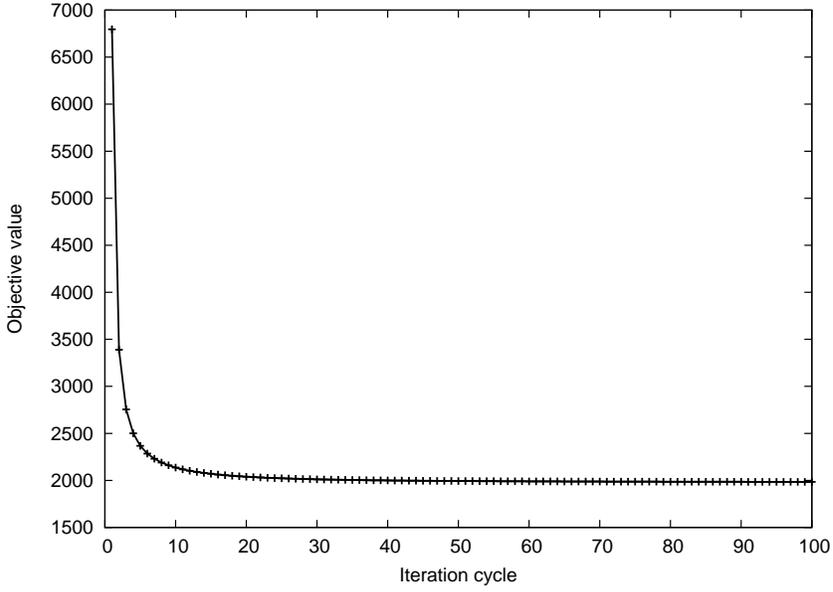}}

\subfloat[Plot of the distance measure (\ref{Normalized  mean absolute distance meas. Eq}).]{\includegraphics[scale=0.9]{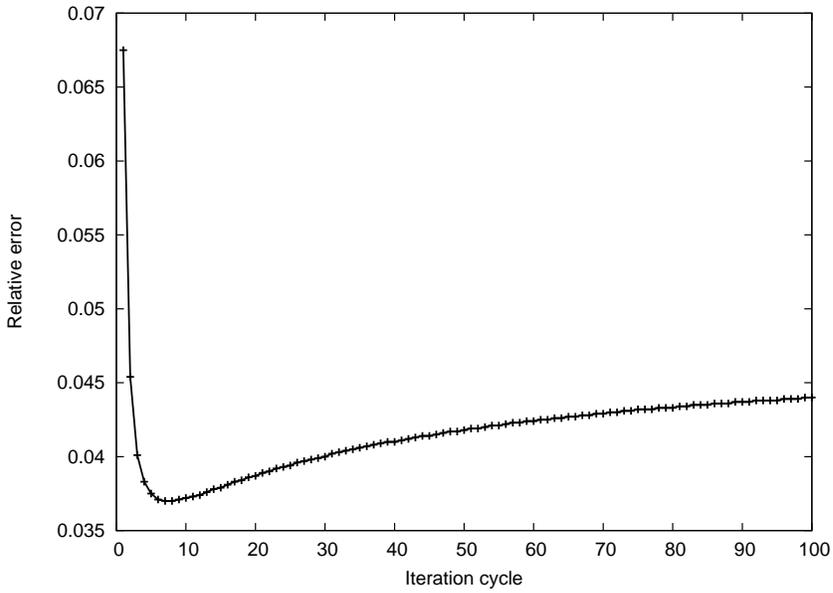}}\caption{\label{fig:Image-reconstruction-by}Image reconstruction by ART when
$M\times N=223,\!744\times51,\!152$.}

\end{figure}

In Section 5.8 of \cite{H09} there is an exact specification of the
so-called \emph{standard projection data} that are used to evaluate
various reconstruction algorithms in that book, the number of lines
used in the standard projection data is $M=223,\!744$. In the evaluations
based on the standard projection data that are reported in \cite{H09}
for reconstruction algorithms that use blob basis functions, the number
of blobs used is $N=51,\!152$. The first experiment on which we report
in this subsection used exactly the same arrangement. (For the experiments
in this subsection, the input data were created and outputs were analyzed
and illustrated using SNARK09 \cite{snark09}.)

In this experiment we applied the ART algorithm of (\ref{U^(0)   , X^(0) is mu_x, U^(k+1),  X^(K+1)  Eq})
and (\ref{c^(k) relate to U^(0)   , X^(0) is mu_x, U^(k+1),  X^(K+1)  Eq})
with $\sigma=5$ and $\lambda=0.05$ to the standard projection data.
In Figure \ref{fig:Image-reconstruction-by}(a) we show the behavior
of the objective function (\ref{eq:minimize}) as a function of iteration
cycles (an \emph{iteration cycle} is defined to be $M$ iterations).
It can be observed that the initial decrease in the objective function
is very rapid.

\begin{figure}
\subfloat[]{\includegraphics[scale=0.35]{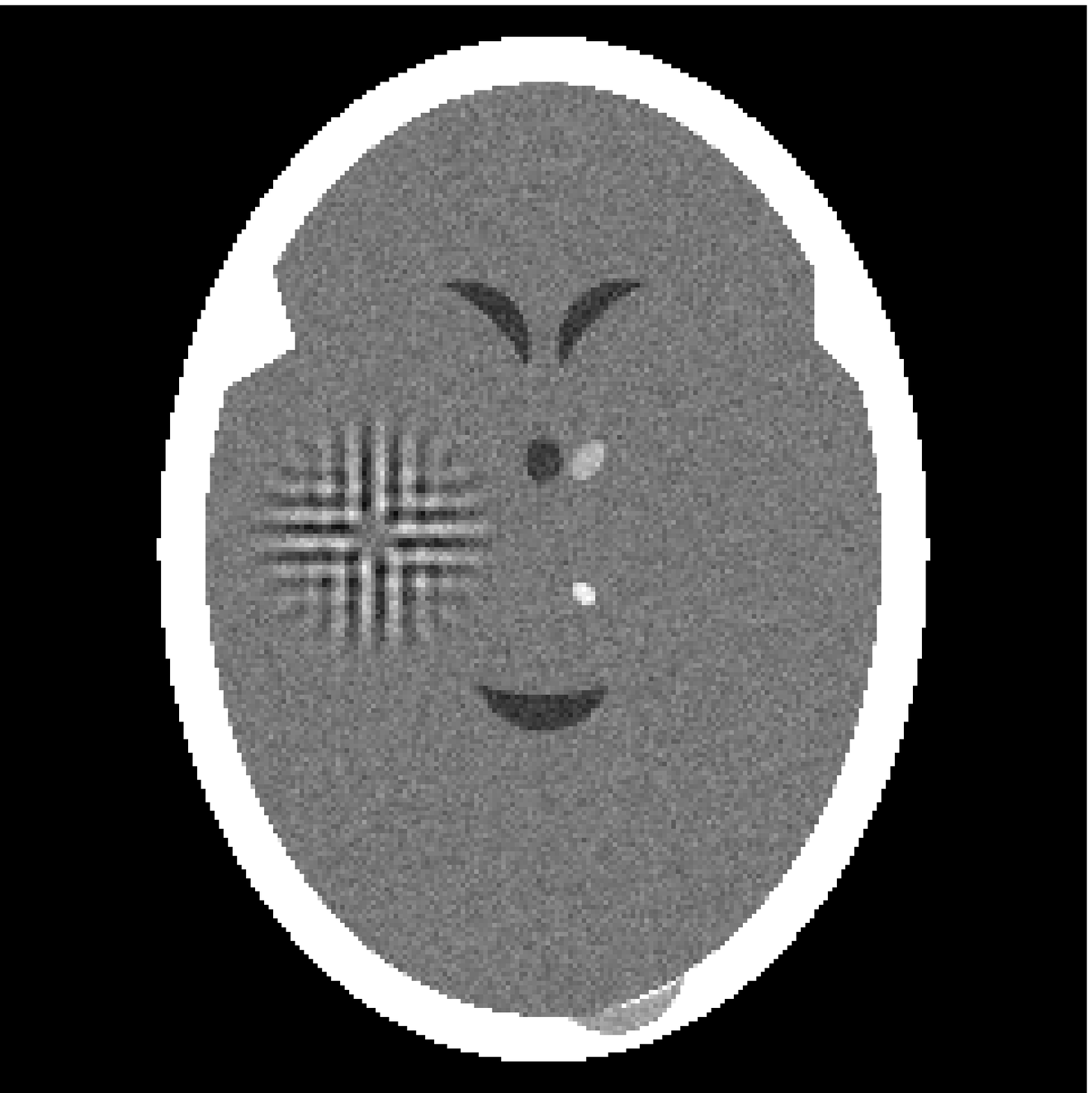}}\subfloat[]{\includegraphics[scale=0.35]{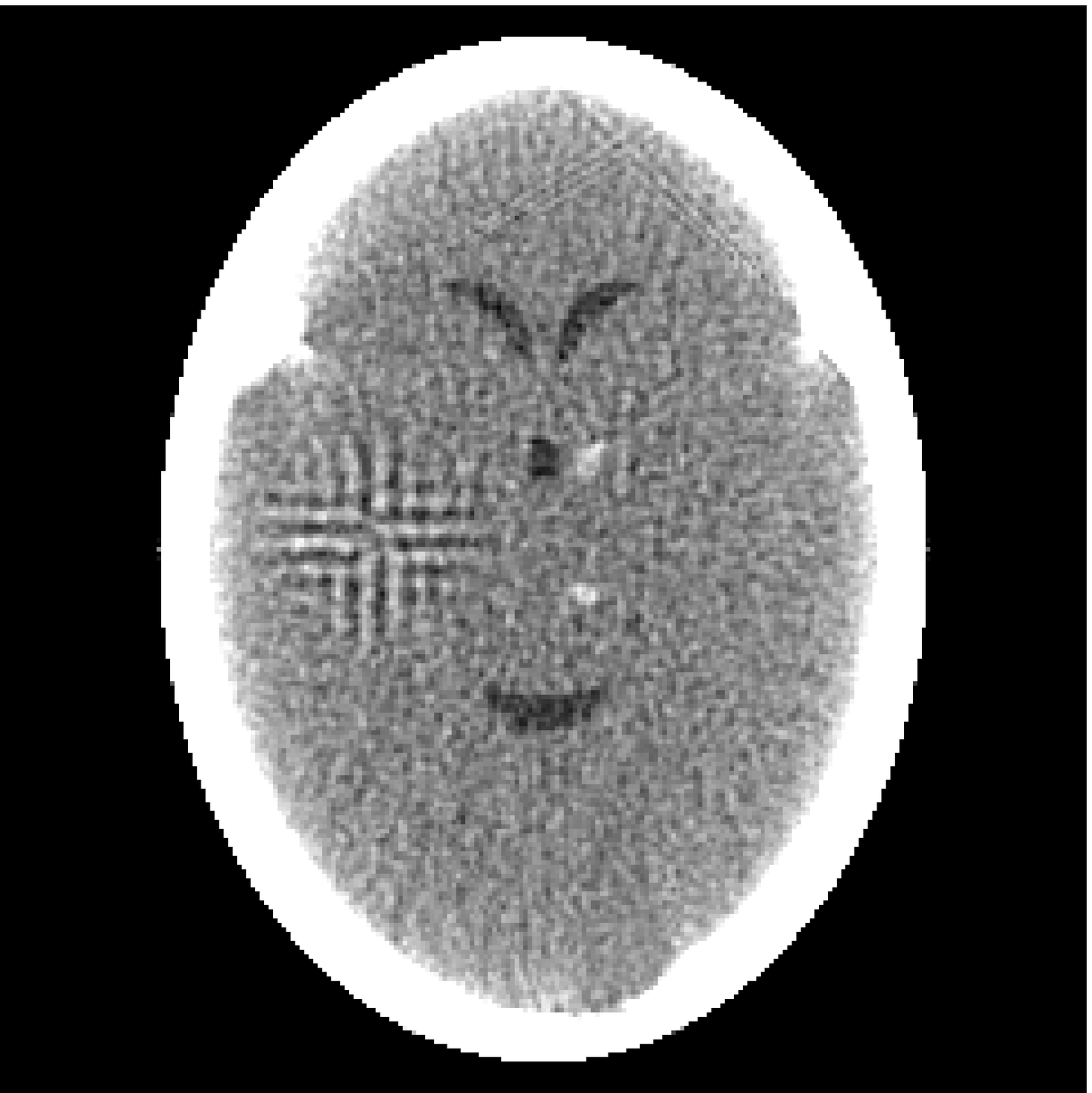}}\caption{\label{fig:Displays-of-}Displays of a $243\times243$ digitized phantom
(a) and of an ART reconstruction when $M\times N=223,\!744\times51,\!152$
(b).}

\end{figure}

This desirable initial behavior is even more noticeable when we evaluate
the algorithm not from the purely mathematical point of view of how
well the objective function is reduced, but rather from the application
point of view of how good are the reconstructed images. For this purpose,
we report on the \emph{normalized mean absolute picture distance measure},
as defined in \cite{H09}. To define this measure we need a $J\times J$
digitization of the test phantom for which the data used in the reconstruction
were collected; such a digitization for the phantom we used is shown
in Figure \ref{fig:Displays-of-}(a). In our definition of the measure
we use $t_{u,v}$ and $s_{u,v}^{(n)}$ to denote the densities of
the \emph{$v$}th pixel of the \emph{$u$}th row of the digitized
test phantom and of the reconstruction (which is obtained from the
vector $x^{(n)}$ of blob coefficients), respectively. We define the
distance measure as\begin{equation}
r^{(n)}=\frac{{\displaystyle \sum_{u=1}^{J}{\displaystyle \sum_{v=1}^{J}\left|t_{u,v}-s_{u,v}^{(n)}\right|}}}{{\displaystyle \sum_{u=1}^{J}{\displaystyle \sum_{v=1}^{J}\left|t_{u,v}\right|}}}.\label{Normalized  mean absolute distance meas. Eq}\end{equation}
In Figure \ref{fig:Image-reconstruction-by}(b) we plot $r^{(n)}$
for this experiment. It is seen that its minimum is reached at the
seventh iteration cycle, i.e., when $n=7M$. This reflects the fact
that the minimization objective (\ref{eq:minimize}) does not (and,
in fact, it cannot in real applications where the phantom is not known
to us) fully describe the application objective. For this reason it
is standard practice in tomography \cite{H09} to stop the iterative
process after a few iteration cycles and use the result at that time
as the reconstruction. The digitization obtained from $x^{(7M)}$
produced by this experiment is shown in Figure \ref{fig:Displays-of-}(b).
The reconstruction is not perfect (as indeed it cannot possibly be
since the measured data are only approximations of the line integrals
assumed by the mathematics), but important features of the phantom
are identifiable in the reconstruction. This ART reconstruction was
carried out in $38.4$ seconds on an Intel Core $1.6$ GHz processor,
$2$ Gbyte memory, $32$ bit laptop.

\begin{figure}
\subfloat[Plot of the objective function (\ref{eq:minimize}).]{\includegraphics[scale=0.9]{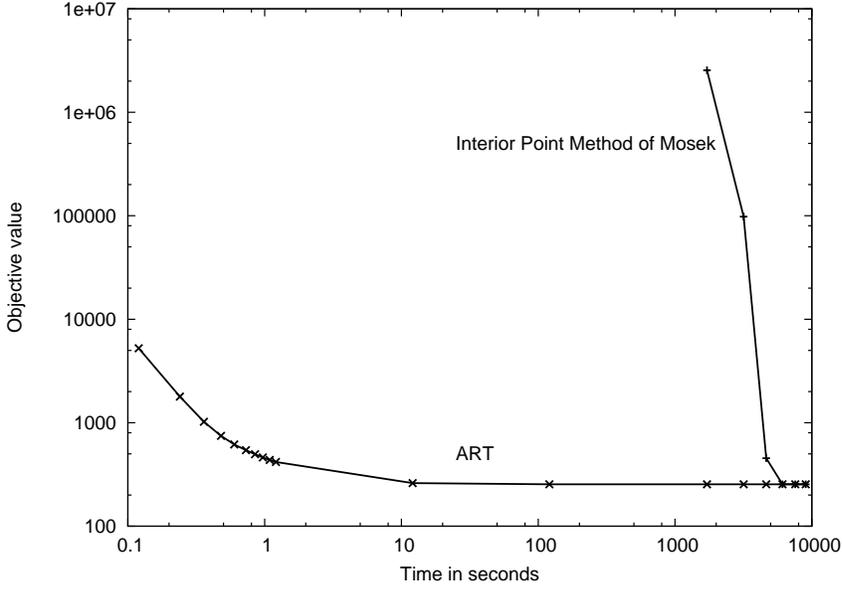}}

$\qquad$\subfloat[Plot of the distance measure (\ref{Normalized  mean absolute distance meas. Eq})]{\includegraphics[scale=0.85]{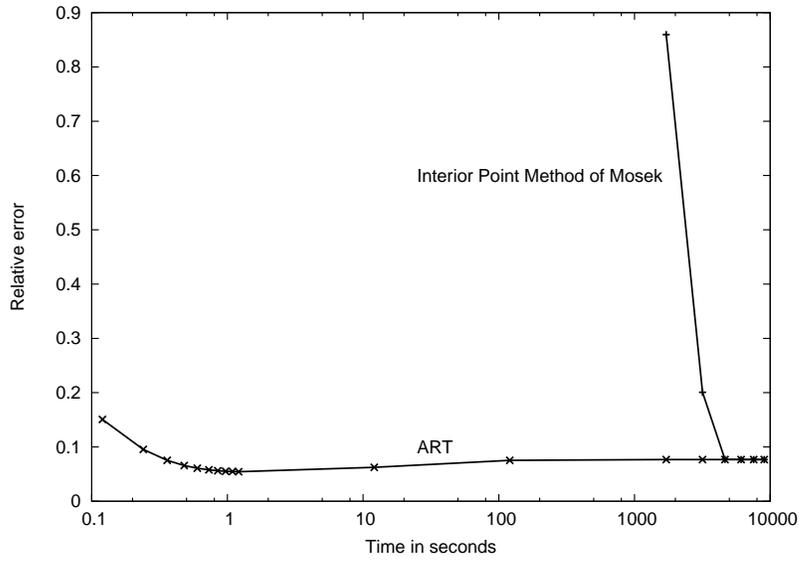}}\caption{\label{fig:Objective}\label{fig:Second_Image-reconstruction}Image
reconstruction by ART and the interior point method in MOSEK when
$M\times N=24,\!880\times5,\!711$.}

\end{figure}

We wanted to compare the time needed by ART with the time needed to
solve the system (\ref{eq:normal}) of consistent equations for the
same data by the current implementation of the interior point method
of MOSEK version 5 \cite{A&A}. Unfortunately this could not be done,
because the memory requirements of the MOSEK software were too large
for our laptop. So we attempted to use a much more powerful Intel
Xeon $2.66$ GHz processor, $16$ Gbyte memory, $64$ bit workstation,
but even the $16$ Gbyte memory was too small to handle this problem
using the MOSEK software. The importance of this memory requirement
issue for the subject matter of this paper cannot be overemphasized:
problems that routinely arise in real applications can be handled
by projection methods using inexpensive laptops, while {}``more sophisticated
alternatives'' fail to produce any results even on much more powerful
workstations due to their much greater demands on computer memory.

In order to be able to compare the efficiency of ART with that of
the interior point method in MOSEK we had to reduce $M$ and $N$
to about a ninth of their previously-used sizes. Thus, in the second
experiment on which we now report $M\times N=24,\!880\times5,\!711$.
For this smaller example we ran both ART and the interior point method
in MOSEK (with its default parameters) on the Intel Xeon $2.66$ GHz
processor, $16$ Gbyte memory, $64$ bit workstation. In Figure \ref{fig:Second_Image-reconstruction}
we plot both the objective function and the distance measure for both
algorithms as a function of time. From the point of view of the objective
function, MOSEK needed over $5000$ seconds to reach a value as low
as ART reached in $10$ seconds. The advantage of ART is more pronounced
when considering the picture distance measure: the optimal value is
reached by ART at $1.7$ seconds (when $n=14M$) while the interior
point method never reaches a distance value that is as low as that
of ART and it needs approximately $5000$ seconds to reach its lowest
distance measure.

\begin{figure}
\subfloat[]{\includegraphics[scale=0.26]{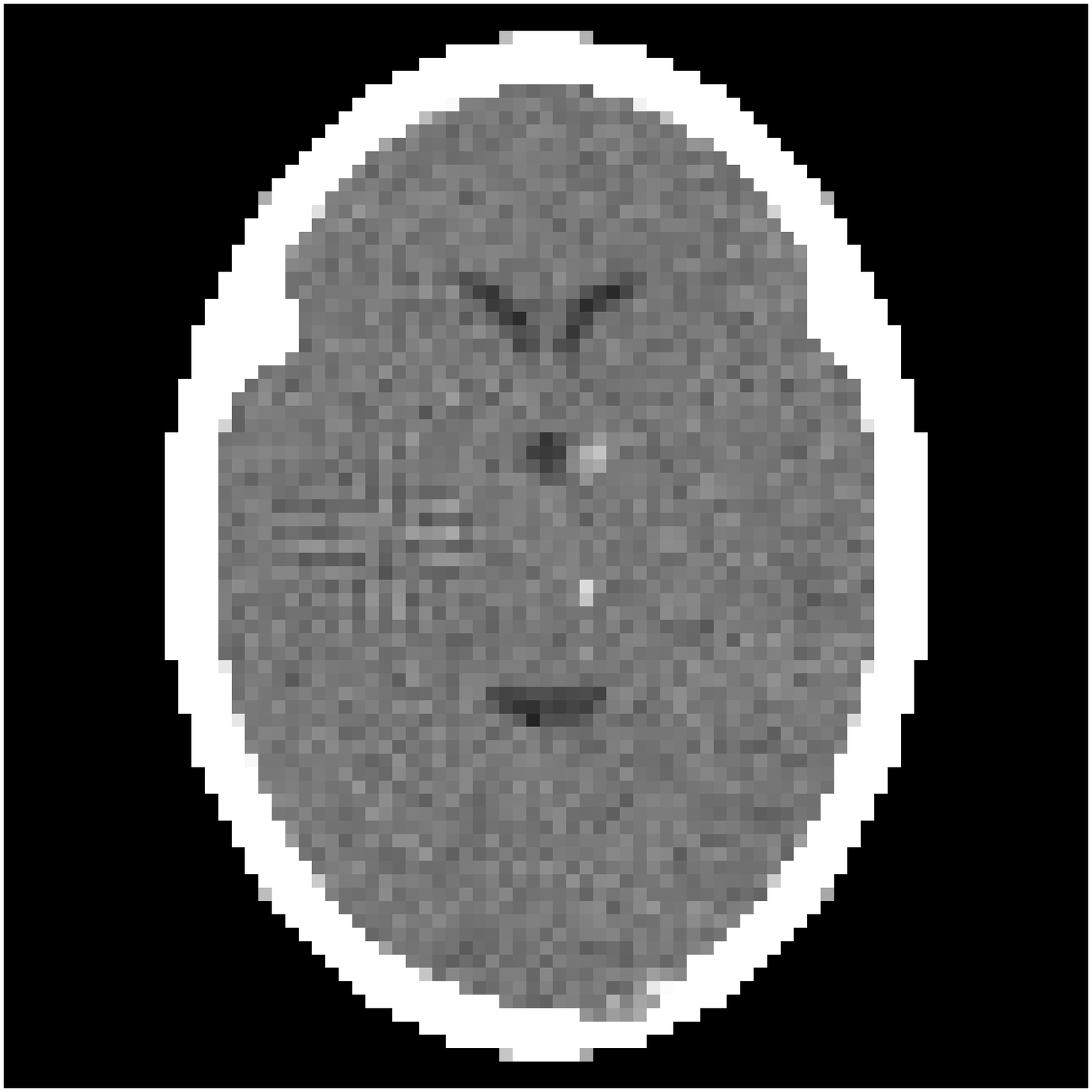}}\subfloat[]{\includegraphics[scale=0.26]{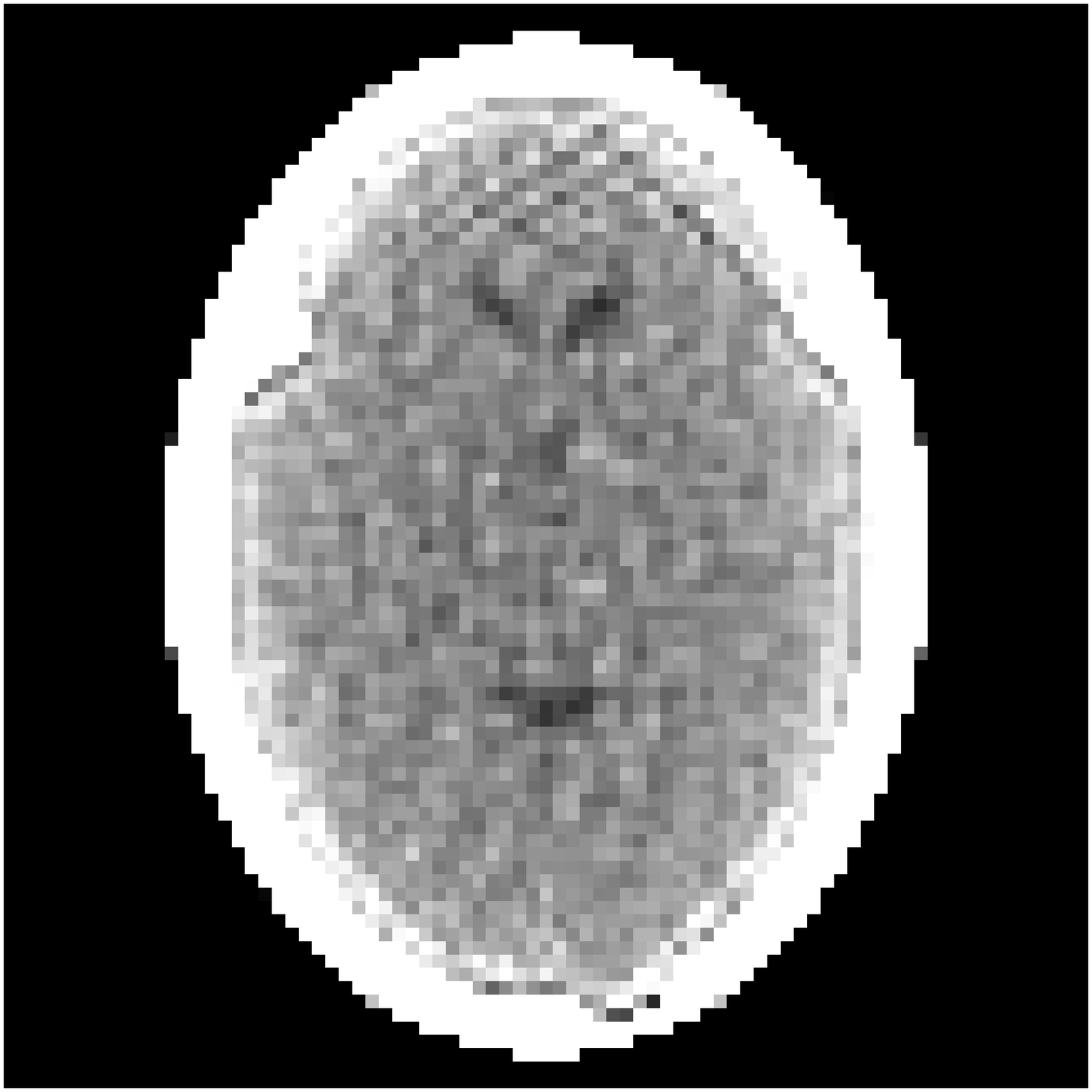}}\caption{\label{fig:Second_Displays-of-}Displays of an $81\times81$ digitized
phantom (a) and of an ART reconstruction when $M\times N=24,\!880\times5,\!711$
(b).}

\end{figure}

Since both $M$ and $N$ are about a ninth of their previous sizes,
we report in Figure \ref{fig:Second_Displays-of-} on the $81\times81$
digitizations of the phantom and of the reconstruction $x^{(14M)}$.
These are clearly inferior to the images in Figure \ref{fig:Displays-of-},
demonstrating the medical necessity for the larger system of equations.

\section{Published and patented results\label{sect:published}}

\subsection{Scientific publications}

Here we give a brief glimpse into some recently published results
that show the efficacy of projection methods for some large problems.
In the problems discussed in \cite{dhc09}, the number of unknowns
was $59,\!049$. In the examples given in \cite{hc08} (a paper devoted
to radiation therapy planning), problems of the form (\ref{eq:blob-representation})
were considered with the number $N$ of unknowns only $515$ but the
number of pairs of constraints $M=128,\!688$. In four out of the
six cases reported there, the projection method ART3+\ \cite{hc08}found
a feasible point in less than three seconds, and in the remaining
two cases a feasible point was found in less than $34$ seconds. These
times are for a standard PC, using an Intel Xeon $1.7$ GHz processor
and $1$ Gbyte memory. The problems in \cite{dhc09,hc08} are small
compared to some of the other applications for which projection methods
have been successfully used. In \cite{chsm06} (a paper devoted to
reconstruction from electron micrographs), there are examples in which
$16,\!777,\!216$ unknowns are to be recovered from $4,\!587,\!520$
measurements (each giving an approximate linear equality) and others
in which $884,\!436$ unknowns are to be recovered from $92,\!160,\!000$
measurements. Projection methods were used in \cite{chsm06} to handle
such large problems in a reasonable time.

In a recent paper \cite{kas09} it is shown that a variant of ART
can be used for crystal lattice orientation distribution function
estimation from diffraction data. One of the problems discussed in
\cite{kas09} has $1,\!372,\!000,\!000$ unknowns and the number of
equations is potentially infinite. They are randomly generated and
a projection step can be carried out as soon as a new equation is
available (an ideal use of a sequential projection method of the row-action
type, see \cite{cen81}). The result reported in the paper for that
problem is the one obtained after $1,\!000,\!000,\!000$ such projection
steps.

As for all methodologies, projection methods are not necessarily the
approach of choice in all applications. However, in important applications
in biomedicine and image processing, projection methods work well
and have been used successfully for a long time. For example, an important
application of reconstruction from projections is electron microscopy
and some of the leading groups in that field consider the projection
method {}``ART with blobs'' to be the method of choice, see \cite{bilbao}.
A mathematical reason for this is that for such problems the angles
between hyperplanes or half-spaces, represented by linear equalities
or linear inequalities as in (\ref{e:feasgould}) and (\ref{eq:blob-representation}),
are in general large (in the sense that the cosine of the angle between
the normals of two randomly chosen hyperplanes in the system to be
solved is likely to be near zero) due to the high sparsity in each
of the rows of the system matrix.

\subsection{Commercial patents}

There is hardly better evidence for the value of projection methods
than the many patents for commercial purposes that include them. Projection
methods are used in commercial devices in many areas. Unfortunately,
if a device is truly commercial, then the algorithm that is actually
used in it is proprietary and usually not published.\ Many commercial
emission tomography scanners use now some sort of iterative algorithms.
A prime example is provided by the commercially-successful Philips
Allegro scanners (see http://www.healthcare.philips.com/main/products/
and \cite{alegro}). In x-ray computerized tomography (CT), there
are reports emanating from companies that sell such scanners indicating
that variants of ART are used in heart imaging; an example is presented
in \cite{isola}.

The first EMI (Electric \& Musical Industries Ltd., London, England,
UK) CT scanner, invented by G.N. Hounsfield \cite{hounsfield}, used
a variant of ART. For this pioneering invention, Hounsfield shared
the Nobel Prize with A.M. Cormack in 1979. Thirty years later (on
September 29, 2009), a patent was issued to Philips (Koninklijke Philips
Electronics N.V., Eindhoven, The Netherlands) for a {}``Method and
device for the iterative reconstruction of cardiac images'' \cite{philips}.
The role of projection methods is demonstrated by the following quote
from the {}``Summary of the Invention'' included in the Patent Description:
\begin{quote}
{}``The iterative reconstruction applied here may particularly be
based on an Algebraic Reconstruction Technique (ART) (cf. R. Gordon,
R. Bender, and G.T. Herman: {}``Algebraic reconstruction techniques
(ART) for three-dimensional electron microscopy and x-ray photography'',
J. Theor. Biol., 29:471--481, 1970) or on a Maximum Likelihood (ML)
algorithm (K. Lange and J.A. Fessler: {}``Globally convergent algorithms
for maximum a posteriori transmission tomography'', IEEE Transactions
on Image Processing, 4(10):1430--1450, 1995), wherein each image update
step uses the projections of a selected subset, i.e., projections
corresponding to a similar movement phase.''
\end{quote}

\section{Conclusion}

In this paper we have shown that, whether or not alternative methods
are applicable, correctly implemented projection methods are very
efficient for convex feasibility problems with linear inequality constraints,
especially for those that are large, sparse, and originate from real-life
applications.\bigskip{}

\textbf{Acknowledgments}. The work of Y. Censor, W. Chen, R. Davidi,
and G.T. Herman is supported by NIH Award Number R01HL070472 from
the National Heart, Lung, and Blood Institute. The content is solely
the responsibility of the authors and does not necessarily represent
the official views of the National Heart, Lung, and Blood Institute
or the National Institutes of Health. The work of P.L. Combettes was
supported by the Agence Nationale de la Recherche under grant ANR-08-BLAN-0294-02.
We thank the Department of Radiation Oncology, Massachusetts General
Hospital and Harvard Medical School (and especially Dr. David Craft)
for providing access to their computers.

\bigskip{}


\begin{thebibliography}{77}
\bibitem{Agmo54}{\small{} Agmon, S.: The relaxation method for linear
inequalities. Canad. J. Math. }\textbf{\small 6}{\small , 382--392
(1954)}{\small \par}

{\small \bibitem{ac89}Aharoni, R., Censor, Y.: Block-iterative projection
methods for parallel computation of solutions to convex feasibility
problems. Linear Algebra Appl. }\textbf{\small 120}{\small , 165--175}\textbf{\small{}
}{\small (1989)}{\small \par}

{\small \bibitem{A&A}Andersen, E.D., Andersen, K.D.: The MOSEK interior
point optimizer for linear programming: an implementation of the homogeneous
algorithm. In: Frenk, H., Roos, K., Terlaky, T., Zhang, S. (eds.)
High Performance Optimization, pp. 197\textendash{}-232, Kluwer, Boston
(2000).}{\small \par}

{\small \bibitem{Andr77} Andrews, H.C., Hunt, B.R.: Digital Image
Restoration. Prentice-Hall, Englewood Cliffs (1977)}{\small \par}

{\small \bibitem{Ausl69} Auslender, A.:}\emph{\small{} }{\small Méthodes
Numériques pour la Résolution des Problèmes d'Optimisation avec Contraintes.
Thèse, Faculté des Sciences, Grenoble (1969)}{\small \par}

\bibitem{Aus76}\textcolor{black}{\small Auslender, A.:}\textcolor{black}{\emph{\small{}
}}\textcolor{black}{\small Optimisation - Méthodes Numériques}\textcolor{black}{\emph{\small .}}\textcolor{black}{\small{}
Masson, Paris (1976)}{\small \par}

\bibitem{b96}{\small{} Bauschke, H.H.: The approximation of fixed
points of compositions of nonexpansive mappings in Hilbert space.
J. Math. Anal. Appl. }\textbf{\small 202}{\small , 150--159 (1996).}{\small \par}

{\small \bibitem{bb96} Bauschke, H.H., Borwein, J.M.: On projection
algorithms for solving convex feasibility problems. SIAM Rev. }\textbf{\small 38}{\small ,
367--426 (1996)}{\small \par}

{\small \bibitem{bc01} Bauschke, H.H., Combettes, P.L.: A weak-to-strong
convergence principle for Fejér-monotone methods in Hilbert spaces.
Math. Oper. Res. }\textbf{\small 26}{\small , 248--264 (2001)}{\small \par}

{\small \bibitem{Baus06} Bauschke, H.H., Combettes, P.L., Kruk, S.G.:
Extrapolation algorithm for affine-convex feasibility problems. Numer.
Algorithms }\textbf{\small 41}{\small , 239--274 (2006)}{\small \par}

\bibitem{Baus03}{\small{} Bauschke, H.H., Deutsch, F., Hundal, H.,
Park, S.-H.: Accelerating the convergence of the method of alternating
projections. Trans. Amer. Math. Soc. }\textbf{\small 355}{\small ,
3433--3461 (2003)}{\small \par}

{\small \bibitem{bmr03} Bauschke, H.H., Matou¨ková, E., Reich, S.:
Projection and proximal point methods: convergence results and counterexamples.
Nonlinear Anal. }\textbf{\small 56}{\small , 715--738 (2004)}{\small \par}

{\small \bibitem{bilbao} Bilbao-Castro, J.R., Marabini, R., Sorzano,
C.O.S., García, I., Carazo, J.M., Fernández, J.J.: Exploiting desktop
supercomputing for three-dimensional electron microscopy reconstructions
using ART with blo}\textcolor{black}{\small bs. J. Struct. Biol. }\textbf{\textcolor{black}{\small 1}}\textbf{\small 65}{\small ,
19--26 (2009)}{\small \par}

{\small \bibitem{Hero06} Blatt, D., Hero, A.O., III: Energy based
sensor network source localization via projection onto convex sets
(POCS). IEEE Trans. Signal Process. }\textbf{\small 54}{\small , 3614--3619
(2006)}{\small \par}

{\small \bibitem{Breg65} Brègman, L.M.: The method of successive
projection for finding a common point of convex sets}\textcolor{black}{\small .
Soviet Math. Dokl. }\textbf{\textcolor{black}{\small 6}}{\small ,
688--692 (1965)}{\small \par}

{\small \bibitem{Butn90} Butnariu, D., Censor, Y.: On the behavior
of a block-iterative projection method for solving convex feasibility
problems. Int. J. Comput. Math. }\textbf{\small 34}{\small , 79--94
(1990)}{\small \par}

\bibitem{Else01}{\small Butnariu, D., Censor, Y., Reich, S. (eds.):
Inherently Parallel Algorithms in Feasibility and Optimization and
Their Applications. Elsevier, Amsterdam (2001)}{\small \par}

\bibitem{bdhk07}{\small Butnariu, D., Davidi, R., Herman, G.T., Kazansev,
I.G.: Stable convergence behavior under summable perturbations of
a class of projection methods for convex feasibility and optimization
problems, IEEE J. Selected Topics Signal Process. }\textbf{\small 1}{\small ,
540--547 (2007)}{\small \par}

{\small \bibitem{chsm06} Carazo, J.M., Herman, G.T., Sorzano, C.O.S.,
Marabini R.: Algorithms for three-dimensional reconstruction from
imperfect projection data provided by electron microscopy. In: Frank,
J. (ed.) }\textit{\emph{\small Electron Tomography: Methods for Three-Dimensional
Visualization of Structures in the Cell}}\emph{\small ,}{\small{} 2nd
edition, pp. 217--243, Springer Science+Business Media, LLC, New York
(2006)}{\small \par}

\bibitem{cen81}{\small Censor, Y.: Row-action methods for huge and
sparse systems and their applications. SIAM Rev. }\textbf{\small 23}{\small ,
444--466 (1981)}{\small \par}

{\small \bibitem{cap88} Censor, Y., Altschuler, M.D., Powlis, W.D.:
On the use of Cimmino's simultaneous projections method for computing
a solution of the inverse problem in radiation therapy treatment planning.
Inverse Problems }\textbf{\small 4}{\small , 607--623 (1988)}{\small \par}

{\small \bibitem{ceh01}Censor, Y., Elfving, T., Herman, G.T.: Averaging
strings of sequential iterations for convex feasibility problems.
In: Butnariu, D., Censor, Y., Reich, S. (eds.) Inherently Parallel
Algorithms in Feasibility and Optimization and Their Applications,
pp. 101--114, Elsevier Science Publishers, Amsterdam (2001)}{\small \par}

\bibitem{cgg01} {\small Censor, Y., Gordon, D., Gordon, R.: BICAV:
A block-iterative, parallel algorithm for sparse systems with pixel-related
weig}\textcolor{black}{\small hting. IEEE Trans. Med. Imaging }\textbf{\textcolor{black}{\small 20}}{\small ,
1050--1060 (2001)}{\small \par}

\bibitem{cs09}{\small Censor, Y., Segal, A.: On the string averaging
method for sparse common fixed points problems. }\textcolor{black}{\small Int.
Trans. Oper. Res. }\textbf{\small 16}{\small , 481--494 (2009)}{\small \par}

{\small \bibitem{ct03}Censor, Y., Tom, E.: Convergence of string-averaging
projection schemes for inconsistent convex feasibility problems. Optim.
Methods Softw. }\textbf{\small 18}{\small , 543--554 (2003)}{\small \par}

{\small \bibitem{CZ97} Censor, Y., Zenios, }\textit{\emph{\small S}}\textit{\small .}\textit{\emph{\small A.:
Parallel Optimization: Theory, Algorithms, and Applications}}\emph{\small .}{\small{}
Oxford University Press, New York (1997)}{\small \par}

{\small \bibitem{Ceti03} Cetin, A.E., Ozaktas, H., Ozaktas, H.M.:
Resolution enhancement of low resolution wavefields with POCS algori}\textcolor{black}{\small thm.
Electron. Lett. }\textbf{\textcolor{black}{\small 3}}\textbf{\small 9}{\small ,
1808--1810 (2003)}{\small \par}

{\small \bibitem{Chen&Herm}Chen, W., Herman, G.T.: Effcient controls
for finitely convergent sequential algorithms. ACM Trans. Math. Software,
to appear}{\small \par}

{\small \bibitem{alegro} Chiang, S., Cardi, C., Matej, S., Zhuang,
H., Newberg, A., Alavi, A., Karp, J.S.: Clinical validation of fully
3-D versus 2.5-D RAMLA reconstruction on the Phillips-ADAC CPET PET
scanne}\textcolor{black}{\small r. Nucl. Med. Commun. }\textbf{\textcolor{black}{\small 25}}\textcolor{black}{\small ,
1103--1107 (2}{\small 004)}{\small \par}

{\small \bibitem{Choi04} Choi, H., Baraniuk, R.G.: Multiple wavelet
basis image denoising using Besov ball projectio}\textcolor{black}{\small ns.
IEEE Signal Process. Lett. }\textbf{\textcolor{black}{\small 11}}\textcolor{black}{\small ,}{\small{}
717--720 (2004)}{\small \par}

\bibitem{Cimm38} {\small Cimmino, G.: Calcolo approssimato per le
soluzioni dei sistemi di equazioni lineari. La Ricerca Scientifica
(Roma) }\textbf{\small 1}{\small , 326--333 (1938)}{\small \par}

{\small \bibitem{Proc93} Combettes, P.L.: The foundations of set
theoretic estimatio}\textcolor{black}{\small n. Proc. IEEE }\textbf{\textcolor{black}{\small 8}}\textbf{\small 1}{\small ,
182--208 (1993)}{\small \par}

{\small \bibitem{c96} Combettes, P.L.: The convex feasibility problem
in image recovery. }\textcolor{black}{\small Adv. Imaging Electron.
Phys. }\textbf{\textcolor{black}{\small 95}}\textcolor{black}{\small ,
}{\small 155--270 (1996)}{\small \par}

{\small \bibitem{Jamo97} Combettes, P.L.: Hilbertian convex feasibility
problem: Convergence of projection methods. Appl. Math. Optim. }\textbf{\small 35}{\small ,
311--330 (1997)}{\small \par}

{\small \bibitem{Imag97} Combettes, P.L.: Convex set theoretic image
recovery by extrapolated iterations of parallel subgradient projections.
IEEE Trans. Image Process. }\textbf{\small 6}{\small , 493--506 (1997)}{\small \par}

\bibitem{Hirs05} {\small Combettes, P.L., Hirstoaga, S.A.: Equilibrium
programming in Hilbert spaces. J. Nonlinear Convex Anal. }\textbf{\small 6}{\small ,
117--136 (2005)}{\small \par}

\bibitem{Sign89} {\small Combettes, P.L., Trussell, H.J.: Methods
for digital restoration of signals degraded by a stochastic impulse
response. IEEE Trans. Acoust. Speech Signal Process. }\textbf{\small 37}{\small ,
393--401 (1989)}{\small \par}

{\small \bibitem{Cott78} Cottle, R.W., Pang, J.-S.: On solving linear
complementarity problems as linear programs.}\textcolor{black}{\small{}
Math. Program. Study }\textbf{\textcolor{black}{\small 7}}\textcolor{black}{\small ,
8}{\small 8--107 (1978)}{\small \par}

\bibitem{crombez02}{\small{} Crombez, G.: Finding common fixed points
of strict paracontractions by averaging strings of sequential iterations,
J. Nonlinear Convex Anal. }\textbf{\small 3,}{\small{} 345--351 (2002)}{\small \par}

{\small \bibitem{dhc09} Davidi, R., Herman, G.T., Censor, Y.: Perturbation-resilient
block-iterative projection methods with application to image reconstruction
from projection}\textcolor{black}{\small s. Int. Trans. Oper. Res.
}\textbf{\small 16}{\small , 505--524 (2009)}{\small \par}

\bibitem{snark09}{\small{} Davidi, R., Herman, G.T., Klukowska, J.:
SNARK09: A programming system for the reconstruction of 2D images
from 1D projections. The CUNY Institute for Software Design and Development.
http://www.snark09.com (2009). Accessed 13 December 2009}{\small \par}

{\small \bibitem{D01} Deutsch, F.: Best Approximation in Inner Product
Spaces. Springer-Verlag, New York (2001)}{\small \par}

\bibitem{Svai09} {\small Eckstein, J., Svaiter, B. F.: General projective
splitting methods for sums of maximal monotone operators. SIAM J.
Control Optim. }\textbf{\small 48}{\small , 787--811 (2009)}{\small \par}

\bibitem{eGG81} {\small Eggermont, P.P.B., Herman, G.T., Lent, A.:
Iterative algorithms for large partitioned linear systems, with applications
to image reconstruction, Linear Algebra Appl. }\textbf{\small 40}{\small ,
37--67, (1981)}{\small \par}

{\small \bibitem{Golu96} Golub, G.H., van Loan, C.F.: Matrix Computations,
3rd ed. Johns Hopkins University Press, Baltimore, MD (1996)}{\small \par}

{\small \bibitem{Gonz01} González-Castaño, F.J., García-Palomares,
U.M., Alba-Castro, J.L., Pousada-Carballo, J.M.: Fast image recovery
using dynamic load balancing in parallel architectures, by means of
incomplete projections. IEEE Trans. Image Process. }\textbf{\small 10}{\small ,
493--499 (2001)}{\small \par}

{\small \bibitem{gbh70}Gordon, R., Bender, R., Herman, G.T.: Algebraic
reconstruction techniques (ART) for three-dimensional electron microscopy
and x-ray photograph}\textcolor{black}{\small y. J. Theor. Biol. }\textbf{\textcolor{black}{\small 2}}\textbf{\small 9}{\small ,
471--482 (1970) }{\small \par}

\bibitem{gould} {\small Gould, N.I.M.: How good are projection methods
for convex feasibility problems? Comput. Optim. Appl. }\textbf{\small 40}{\small ,
1--12 (2008)}{\small \par}

{\small \bibitem{Gu04} Gu, J., Stark, H., Yang, Y.: Wide-band smart
antenna design using vector space projection methods. IEEE Trans.
Antennas Propag. }\textbf{\small 52}{\small , 3228--3236 (2004)}{\small \par}

{\small \bibitem{Gubi67} Gubin, L.G., Polyak, B.T., Raik, E.V.: The
method of projections for finding the common point of convex sets.
USSR Com}\textcolor{black}{\small put. Math. Math. Phys. }\textbf{\textcolor{black}{\small 7}}\textcolor{black}{\small ,
1-}{\small -24 (1967)}{\small \par}

{\small \bibitem{H09} Herman, G.T.:}\emph{\small{} }\textit{\emph{\small Fundamentals
of Computerized Tomography: Image Reconstruction from Projections}}\emph{\small ,
}{\small 2nd ed. Springer, New York (2009)}{\small \par}

{\small \bibitem{hc08} Herman, G.T., Chen, W.: A fast algorithm for
solving a linear feasibility problem with application to intensity-modulated
radiation therapy. Linear Algebra Appl. }\textbf{\small 428}{\small ,
1207--1217}\textbf{\small{} }{\small (2008)}{\small \par}

{\small \bibitem{hounsfield} Hounsfield, G.N.: A method and apparatus
for examination of a body by radiation such as X or gamma radiation.
UK Patent No. 1283915 (1968/72)}{\small \par}

{\small \bibitem{isola} Isola, A.A., Ziegler, A., Koehler, T., Niessen,
W.J., Grass, M.: Motion compensated iterative cone-beam CT image reconstruction
with adapted blobs as basis functio}\textcolor{black}{\small ns. Phys.
Med. Biol. }\textbf{\textcolor{black}{\small 5}}\textbf{\small 3}{\small ,
6777--6797 (2008)}{\small \par}

\bibitem{Kacz37} {\small Kaczmarz, S.: Angen\"{a}herte Aufl\"{o}sung
von Systemen linearer Gleichungen. Bull. Acad. Sci. Pologne }\textbf{\small A35}{\small ,
355--357 (1937)}{\small \par}

{\small \bibitem{kas09} Kazantsev, I.G., Schmidt, S., Poulsen, H.F.:
A discrete spherical x-ray transform of orientation distribution functions
using bounding cubes. Inverse Problems }\textbf{\small 25}{\small ,
105009 (2009)}{\small \par}

{\small \bibitem{Lopu97} Kiwiel, K.C., {\L{}}opuch, B.: Surrogate
projection methods for finding fixed points of firmly nonexpansive
mappings. SIAM J. Optim. }\textbf{\small 7}{\small , 1084--1102 (1997)}{\small \par}

{\small \bibitem{Lee08} Lee, S.-H., Kwon, K.-R.: Mesh watermarking
based on projection onto two convex sets. }\textcolor{black}{\small Multimedia
Syst. }\textbf{\textcolor{black}{\small 13}}\textcolor{black}{\small ,
323}{\small --330 (2008)}{\small \par}

{\small \bibitem{Lew90} Lewitt, R.M.: Multidimensional digital image
representation using generalized Kaiser-Bessel window functions. J.
Opt. Soc. Amer. A }\textbf{\small 7}{\small , 1834--1846 (1990)}{\small \par}

\bibitem{Liew05}{\small Liew, A.W.-C., Yan, H., Law, N.-F.: POCS-based
blocking artifacts suppression using a smoothness constraint set with
explicit region modeling}\textcolor{black}{\small . IEEE Trans. Circuits
Syst. Video Technol. }\textbf{\textcolor{black}{\small 15}}\textcolor{black}{\small ,
795--800 (2005)}{\small \par}

\bibitem{Lu09}\textcolor{black}{\small Lu, Y.M., Karzand, M., Vetterli,
M.: Demosaickin}{\small g by alternating projections: Theory and fast
one-step implementation. IEEE Trans. Image Process., to appear}{\small \par}

{\small \bibitem{Merz63} Merzlyakov, Y.I.: On a relaxation method
of solving systems of linear inequalities. }\textcolor{black}{\small USSR
Comput. Math. Math. Phys. }\textbf{\textcolor{black}{\small 2}}\textcolor{black}{\small ,
504--510 (1963)}{\small \par}

\bibitem{Motz54} {\small Motzkin, T.S., Schoenberg, I.J.: The relaxation
method for linear inequalities. Canad. J. Math. }\textbf{\small 6}{\small ,
393--404 (1954)}{\small \par}

{\small \bibitem{Otta88} Ottavy, N.: Strong convergence of projection-like
methods in Hilbert spaces. J. Optim. Theory Appl. }\textbf{\small 56}{\small ,
433--461 (1988)}{\small \par}

\bibitem{pen09} {\small Penfold, S.N., Schulte, R.W., Censor, Y.,
Bashkirov, V., McAllister, S., Schubert, K.E., Rosenfeld, A.B.: Block-iterative
and string-averaging projection algorithms in proton computed tomography
image reconstruction. In: Censor, Y., Jiang, M., Wang, G. (eds.) Biomedical
Mathematics: Promising Directions in Imaging, Therapy Planning and
Inverse Problems, Medical Physics Publishing, Madison, WI, to appear.}{\small \par}

{\small \bibitem{Pier76} Pierra, G.: Éclatement de contraintes en
parallèle pour la minimisation d'une forme quadratique. In: Proc.
7th IFIP Conf. on Optimization Techniques: Modeling and Optimization
in the Service of Man, Lecture Notes in Comput. Sci. }\textbf{\small 41}{\small ,
pp. 200--218. Springer-Verlag, London (1976)}{\small \par}

\bibitem{Pier84}{\small Pierra, G.: Decomposition through formalization
in a product space. Math. Programming }\textbf{\small 28}{\small ,
96--115 (1984)}{\small \par}

\bibitem{rh03}{\small Rhee, H.: An application of the string averaging
method to one-sided best simultaneous approximation. J. Korea Soc.
Math. Educ. Ser. B: Pure Appl. Math. }\textbf{\small 10}{\small ,
49--56 (2003)}{\small \par}

{\small \bibitem{Sams04} Samsonov, A.A., Kholmovski, E.G., Parker,
D.L., Johnson, C.R.: POCSENSE: POCS-based reconstruction for sensitivity
encoded magnetic resonance imagin}\textcolor{black}{\small g. Magn.
Reson. Med. }\textbf{\textcolor{black}{\small 52}}\textcolor{black}{\small ,
1397--1406 (2004)}{\small \par}

{\small \bibitem{Shak08} Shaked, N.T., Rosen, J.: Multiple-viewpoint
projection holograms synthesized by spatially incoherent correlation
with broadband functions. J. Opt. Soc. Amer. A }\textbf{\small 25}{\small ,
2129--2138 (2008)}{\small \par}

{\small \bibitem{Shar00} Sharma, G.: Set theoretic estimation for
problems in subtractive color}\textcolor{black}{\small . Color Res.
Appl.}\textcolor{red}{\small{} }\textbf{\small 25}{\small , 333-348
(2000)}{\small \par}

{\small \bibitem{Star98} Stark, H., Yang, Y.: Vector Space Projections
: A Numerical Approach to Signal and Image Processing, Neural Nets,
and Optics. Wiley-Interscience, New York (1998)}{\small \par}

{\small \bibitem{Vwik04}van Wyk, B.J., van Wyk, M.A.: A POCS-based
graph matching algorith}\textcolor{black}{\small m. IEEE Trans. Pattern
Anal. Machine Intell. }\textbf{\textcolor{black}{\small 26}}\textcolor{black}{\small ,
1526--1530 (2004)}{\small \par}

{\small \bibitem{Youl82} Youla, D.C., Webb, H.: Image restoration
by the method of convex projections: Part 1 -- theo}\textcolor{black}{\small ry.
IEEE Trans. Med. Imaging }\textbf{\textcolor{black}{\small 1}}\textcolor{black}{\small ,
81--94 }{\small (1982)}{\small \par}

{\small \bibitem{Yuka06} Yukawa, M., Yamada, I.: Pairwise optimal
weight realization -- Acceleration technique for set-theoretic adaptive
parallel subgradient projection algorithm. IEEE Trans. Signal Process.
}\textbf{\small 54}{\small , 4557--4571 (2006)}{\small \par}

\bibitem{Zhan05}{\small Zhang, T., Hong, H.: Restoration algorithms
for turbulence-degraded images based on optimized estimation of discrete
values of overall point spread functions. Opt. Eng. }\textbf{\small 44}{\small ,
017005 (2005)}{\small \par}

{\small \bibitem{philips} Ziegler, A., Grass, M., Koehler, T.: Method
and device for the iterative reconstruction of cardiac images. US
Patent Number 7596204 (2009)}
\end{thebibliography}
\end{document}